\documentclass{amsart}
 \usepackage{amssymb,amsmath,amsfonts,epsfig,latexsym}

 \newtheorem{theorem}{Theorem}[section]
\newtheorem{definition}[theorem]{Definition}
\newtheorem{proposition}[theorem]{Proposition}
\newtheorem{lemma}[theorem]{Lemma}
\newtheorem{corollary}[theorem]{Corollary}

\newtheorem*{BF}{Brion's Formula}

\theoremstyle{definition}
\newtheorem{remark}[theorem]{Remark}
\newtheorem{example}[theorem]{Example}

\def\Z{\ensuremath{\mathbb{Z}}}
\def\Q{\ensuremath{\mathbb{Q}}}
\def\P{\ensuremath{\mathbb{P}}}

\def\R{\ensuremath{\mathbb{R}}}

\def\E{\ensuremath{\mathcal{E}}}
\def\O{\ensuremath{\mathcal{O}}}

\def\m{\ensuremath{\mathfrak{m}}}

\def\bu{\ensuremath{\mathbf{u}}}

\def\<{\ensuremath{\langle}}
\def\>{\ensuremath{\rangle}}

\def\PP{P\!P}

\DeclareMathOperator{\divisor}{div}
\DeclareMathOperator{\Hilb}{Hilb}
\DeclareMathOperator{\Hom}{Hom}

\DeclareMathOperator{\Pic}{Pic}
\DeclareMathOperator{\rk}{rk}
\DeclareMathOperator{\Span}{span}
\DeclareMathOperator{\Spec}{Spec}
\DeclareMathOperator{\Star}{Star}
\DeclareMathOperator{\Sym}{Sym}
\DeclareMathOperator{\Poly}{Poly}

\begin{document}

\title[Piecewise polynomials and Minkowski weights]{Piecewise polynomials, {M}inkowski weights, and localization on toric varieties}

\author[Katz]{Eric Katz}
\address{ Department of Mathematics, 1 University Station C1200, Austin, TX 78712}
\email{eekatz@math.utexas.edu}
\author[Payne]{Sam Payne$^*$}\thanks{$^*$Supported by the Clay Mathematics Institute.  Part of this research was done during a visit to the Institut Mittag-Leffler (Djursholm, Sweden).}
\address{Stanford University, Mathematics, Bldg. 380, 450 Serra Mall, Stanford, CA 94305}
\email{spayne@stanford.edu}

\begin{abstract}
We use localization to describe the restriction map from equivariant Chow cohomology to ordinary Chow cohomology for complete toric varieties in terms of piecewise polynomial functions and Minkowski weights.  We compute examples showing that this map is not surjective in general, and that its kernel is not always generated in degree one.  We prove a localization formula for mixed volumes of lattice polytopes and, more generally, a Bott residue formula for toric vector bundles.
\end{abstract}

\maketitle

\tableofcontents

\section{Introduction}

Let $\Delta$ be a complete fan in $N_\R$, where $N$ is a lattice of rank $n$, and let $X = X(\Delta)$ be the corresponding complete $n$-dimensional toric variety.  See \cite{Fulton93} for standard notation and general background on toric varieties.  The equivariant operational Chow cohomology ring with integer coefficients $A^*_T(X)$ is naturally isomorphic to the ring of integral piecewise polynomial functions on $\Delta$ \cite{chowcohom}, and there is a canonical map to ordinary Chow cohomology with integer coefficients
\[
\iota^*: A^*_T(X) \rightarrow A^*(X)
\]
induced by inclusions of $X$ in the finite dimensional approximations of the Borel mixed space \cite{EdidinGraham98}.  Now $A^*(X)$ is naturally isomorphic to the ring of Minkowski weights on $\Delta$ \cite{FultonSturmfels97}, and $\iota^*$ has a natural interpretation in terms of localization and equivariant multiplicities, as follows.

Let $M = \Hom(N, \Z)$, which is naturally identified with the character lattice of $T$, and let $\Sym^\pm(M)$ be
the $\Z$-graded ring obtained by inverting all of the homogeneous
elements in the ring $\Sym^*(M)$ of polynomials with integer coefficients.  We refer to
elements of $\Sym^\pm(M)$ as rational functions, and elements of
the subring $\Sym^*(M)$ as polynomials.  Each maximal cone $\sigma
\in \Delta$ corresponds to a nondegenerate torus fixed point
$x_\sigma \in X$,  which has an ``equivariant multiplicity"
$e_{x_\sigma}[X] \in \Sym^\pm (M)$, which is a homogeneous
rational function of degree $- n$.  Since every rational
polyhedral cone admits a unimodular subdivision, these equivariant
multiplicities are determined by the following two properties.
\begin{enumerate}
\item  If $\sigma_1, \ldots, \sigma_r$ are the maximal cones of a rational polyhedral subdivision of a cone $\sigma$, then
\[
e_{x_\sigma}[X] = e_{x_{\sigma_1}}[X] + \cdots + e_{x_{\sigma_r}}[X].
\]
\item  If $\sigma$ is a unimodular cone, spanned by a basis $e_1, \ldots, e_n$ for $N$, then
\[
e_{x_\sigma}[X] = \frac{1}{e_1^* \cdots e_n^*} \ .
\]
\end{enumerate}
The fact that the sum of rational functions determined by (1) and
(2) is independent of the choice of unimodular subdivision is not
obvious from elementary considerations, though it follows directly from the theory of localization at torus fixed
points in algebraic geometry \cite{EdidinGraham98b} and the theory of
equvariant multiplicities developed by Rossmann \cite{Rossmann89} and
Brion.  See, in particular, Theorem~4.2 and Proposition~4.3 of
\cite{Brion97}. Here we give a
combinatorial proof of this independence; the techniques of this
proof may be of independent interest.  We view the multigraded
Hilbert function $\Hilb(\sigma)$ of the affine toric variety $U_\sigma$, given by,
\[
\Hilb(\sigma) = \sum_{u \in (\sigma^* \cap M)} x^u
\]
as a rational function on the dense torus $T \subset X$.  We define
$e_\sigma$ to be $(-1)^n$ times the quotient of the leading forms
when $\Hilb(\sigma)$ is written as a quotient of two polynomials
in local coordinates at the identity $1_T$.  We then show that
$e_\sigma$ satisfies properties analogous to (1) and (2) and
therefore is equal to $e_{x_\sigma}[X]$.  See Section
\ref{combinatorics} for details. Our
approach is inspired by the presentation of multidegrees
of multigraded modules over polynomial rings in \cite[Sections~1.2 and 1.7]{KnutsonMiller05} and \cite[Chapter
8]{MillerSturmfels05}.

Recall that the ring of integral piecewise polynomial functions
$\PP^*(\Delta)$ is the ring of continuous functions $f: |\Delta|
\rightarrow \R$ such that the restriction $f_\sigma$ of $f$ to
each maximal cone $\sigma \in \Delta$ is a polynomial in
$\Sym^*(M)$.

\begin{proposition} \label{pushforward}
Let $\Delta$ be a complete $n$-dimensional fan, and let $f \in \PP^k(\Delta)$ be a piecewise polynomial function.  Then
\[
\sum_{\sigma} e_{\sigma} f_\sigma
\]
is a homogeneous polynomial in $\Sym^*(M)$ of degree $k - n$.
\end{proposition}

\noindent In particular, if the degree of $f$ is less than $n$ then $\sum e_\sigma f_\sigma$ vanishes.  If $\deg f = n$, then $\sum e_\sigma f_\sigma = d$ is an integer, which may be identified with the codimension $n$ Minkowski weight $c(0) = d$ on $\Delta$.

Minkowski weights of codimension less than $n$ may be constructed
similarly from piecewise polynomials using equivariant
multiplicities, as follows.  For any cone $\tau \in \Delta$, let
$\Delta_\tau$ be the fan in $\big(N / (N \cap \Span \tau)\big)_\R$
whose cones are the projections of the cones in $\Delta$ that
contain $\tau$.  If $\sigma$ is a maximal such cone, we define
$e_{\sigma, \tau}$ to be $e_{\overline{\sigma}}$, where
$\overline{\sigma}$ is the image of $\sigma$ in $\Delta_\tau$.  So
$e_{\sigma, \tau}$  is a homogeneous rational function of degree
$(\dim \tau - n)$ in the graded subring $\Sym^\pm(\tau^\perp \cap
M)$ of $\Sym^\pm (M)$.

\begin{proposition} \label{restricted pushforward}
Let $\Delta$ be a complete fan, and let $f \in \PP^k(\Delta)$ be a piecewise polynomial function.  Then, for any $\tau \in \Delta$,
\[
c(\tau) = \sum_{\sigma \succeq \tau} e_{\sigma, \tau} f_\sigma
\]
is a homogeneous polynomial in $\Sym^*(M)$ of degree $k  + \dim \tau - n$.
\end{proposition}

\noindent If $k \leq n$ then $c(\tau)$ is an integer for every codimension~$k$ cone in $\Delta$, and these integers are a Minkowski weight of codimension~$k$.  Propositions~\ref{pushforward} and \ref{restricted pushforward} are proved in Section~\ref{localizationsection} using elementary properties of generating functions for lattice points in polyhedral cones.

 \begin{remark}
Proposition~\ref{pushforward} is the special case of
Proposition~\ref{restricted pushforward} where $\tau = 0$.  The
essential content of Propositions~\ref{pushforward} and
\ref{restricted pushforward} is that the denominator of the sum
must divide the numerator.  In some special cases, this
divisibility may be seen as a consequence of Brion's Formula, and
its generalizations, in the theory of generating functions for
lattice points in polyhedra. See Section \ref{mixed volumes}
below.  Other special cases of these cancellations appeared
earlier in \cite{Brion96}; in particular, Brion showed that $\sum
e_{\sigma, \tau} f_\sigma$ is in $\Sym^*(M_\Q)$ when $\Delta$ is
simplicial.
 \end{remark}

 \begin{theorem} \label{localization}
 The natural map $\iota^*: A^k_T(X) \rightarrow A^k(X)$ takes the equivariant Chow cohomology class corresponding to a piecewise polynomial function $f$ to the ordinary Chow cohomology class corresponding to the Minkowski weight $c$ given by
 \[
 c(\tau) = \sum_{\sigma \succeq \tau} e_{\sigma, \tau} f_\sigma,
 \]
 for all codimension $k$ cones $\tau \in \Delta$.
\end{theorem}

\noindent We prove Theorem~\ref{localization} in Section~\ref{localizationsection} by interpreting Propositions~\ref{pushforward} and \ref{restricted pushforward} in terms of general localization formulas in equivariant Chow cohomology \cite{EdidinGraham98b}.

We apply Theorem~\ref{localization} to study the map $\iota^*: A^*_T(X) \rightarrow A^*(X)$.  Recall that if $X$ is smooth then capping with the fundamental class of $X$ gives isomorphisms
\[
A^*(X) \cong A_{n-*}(X) \mbox{ \ and \ } A^*_T(X) \cong A_{n-*}^T(X).
\]
Furthermore, the globally linear functions $u \in M$, identified with the equivariant first Chern classes of the toric line bundles $\O (\divisor \chi^{u})$, act on $A_{*}^T(X)$ as homogeneous operators of degree $-1$, and there is a natural isomorphism to ordinary Chow homology \cite[Section~2.3]{Brion97},
\[
A_*^T(X) / M A_*^T(X) \cong A_*(X).
\]
It follows that if $X$ is a smooth toric variety then $\iota^*$ is surjective and its kernel is generated by $M$ in degree one.  Similar arguments show that if $X$ is simplicial, then $\iota^*$ becomes surjective after tensoring with $\Q$, with kernel generated by $M_\Q$ in degree one.

\begin{theorem} \label{surface}
There exist projective toric surfaces $X$ such that $\iota^*: A^2_T(X) \rightarrow A^2(X)$ is not surjective.
\end{theorem}

\noindent In particular, even when the natural map $\PP^*(\Delta)/(M) \rightarrow A^*(X)$ becomes an isomorphism after tensoring with $\Q$, it need not be an isomorphism over $\Z$.

\begin{theorem} \label{threefold}
There exist projective toric threefolds $X$ such that $\iota^*:A^*_T(X)_\Q \rightarrow A^*(X)_\Q$ is not surjective and its kernel is not generated in degree one.
\end{theorem}

\noindent It follows that the natural map from piecewise polynomials modulo linear functions to Minkowski weights is neither injective nor surjective in general.  We prove Theorems~\ref{surface} and \ref{threefold} in Section~\ref{applications} by computing the maps $A^*_T(X) \rightarrow A^*(X)$ for several examples of singular toric varities using Theorem~\ref{localization}.

\begin{remark}
Minkowski weights on $\Delta$, and classes in $A^*(X)$, correspond to tropical varieties supported on the cones of $\Delta$, and are of significant interest in tropical geometry \cite[Section~9]{Katz06} \cite[p.~10]{Mikhalkin06}.  The desire to use piecewise polynomials to produce interesting examples of Minkowski weights was one of the main motivations for this research.  We hope and expect that the combinatorial localization techniques developed here will be useful in tropical geometry.
\end{remark}

\noindent \textbf{Acknowledgments.}  We thank B.~Sturmfels for illuminating the connection between Minkowksi weights and tropical geometry, D.~Speyer for enlightening discussions on Brion's Formula, E.~Miller for helping us to understand multidegrees, and T.~Ekedahl for useful conversations related to principal parts of rational functions.

\section{Combinatorics of equivariant multiplicities} \label{combinatorics}

Let $N$ be a lattice of rank $n$, and let $M = \Hom(N,\Z)$ be its dual lattice.  Let $\Poly(N)$ denote the rational polytope algebra on $N_\R$, the subring of real-valued functions on $N_\R$ generated by the characteristic functions of closed rational polyhedra.  We write $[Q] \in \Poly(N)$ for the characteristic function of a closed polyhderon $Q$.  Recall that $Q$ has a polar dual $Q^*$, which is a closed polyhedron in $M_\R$, defined by
\[
Q^* = \{ u \in M_\R \, | \, \< u,v \> \geq -1 \mbox{ for all } v \in Q \},
\]
and there is a linear map from $\Poly(N)$ to $\Poly(M)$ given by $[Q] \mapsto [Q^*]$  \cite{Lawrence88}.  Furthermore, there is a linear map
\[
\nu: \Poly(M) \rightarrow \Q(M),
\]
to the quotient field $\Q(M)$ of the multivariate Laurent polynomial ring $\Z[M]$, that takes the class of a closed, pointed polyhedron $P$ to the generating function $\sum_{u \in (P \cap M)} x^u$, expressed as a rational function, and takes the class of a polyhedron containing a line to 0 \cite[Theorem~VIII.3.3]{Barvinok02}.  In particular, for any closed polyhedral cone $\sigma$ in $N_\R$, $\nu(\sigma^*) = \Hilb(\sigma)$, where
\[
\Hilb(\sigma) = \sum_{u \in (\sigma^* \cap M)} x^u
\]
is the multigraded Hilbert series of the affine toric variety $U_\sigma$.  Composing polar duality with the valuation $\nu$ then gives a linear map
\[
\nu^*: \Poly(N) \rightarrow \Q(M)
\]
that takes $[\sigma]$ to $\Hilb(\sigma)$.

\begin{lemma} \label{additive hilb}
If $\sigma_1, \ldots, \sigma_r$ are the maximal cones in a rational polyhedral subdivision of an $n$-dimensional cone $\sigma$, then
\[
\Hilb(\sigma) = \Hilb(\sigma_1) + \cdots + \Hilb(\sigma_r).
\]
\end{lemma}

\begin{proof}
In the polytope algebra $\Poly(N)$,
\[
[\sigma] = [\sigma_1] + \cdots + [\sigma_r] \pm  \mbox{ classes of lower dimensional cones. }
\]
Since the duals of lower dimensional cones contain lines, these terms are all in the kernel of $\nu^*$.  Therefore, $\nu^*([\sigma]) = \nu^*([\sigma_1]) + \cdots + \nu^*([\sigma_r])$, and the lemma follows.
\end{proof}

The generating function $\Hilb(\sigma)$, being an element of $\Q(M)$, is naturally interpreted as a rational function on the torus $T = \Spec \Q[M]$.  Therefore, $\Hilb(\sigma)$ may be expanded as a quotient of two power series in local parameters at the identity $1_T$.  The principal part of this expansion, the quotient of the leading forms, which we denote by
\[
\Hilb(\sigma)_\circ \in \Sym^\pm(M_\Q),
\]
is a rational function on the tangent space of $T$ at $1_T$, which cuts out the tangent cone of zeros of $\Hilb(\sigma)$ minus the tangent cone of its poles.  See the Appendix for details on principal parts of rational functions.

\begin{definition} \label{multiplicity}
If  $\sigma$ is an $n$-dimensional rational polyhedral cone in $N_\R$ then
\[
e_\sigma = (-1)^n \cdot \Hilb(\sigma)_\circ.
\]
\end{definition}

\begin{lemma} \label{degree hilb}
If $\sigma$ is an $n$-dimensional rational polyhedral cone in $N_\R$, then $e_\sigma$ is homogeneous of degree $-n$.
\end{lemma}

\begin{proof}
The lemma follows directly from closed formulas for polyhedral generating functions, such as those given in \cite{ehrhartseries}, as follows.  Suppose $\Sigma$ is a unimodular subdivision of $\sigma^*$, and $u_1, \ldots, u_s$ are the primitive generators of the rays of $\Sigma$.  Then every lattice point in $\sigma^*$ lies in the relative interior of a unique cone $\tau \in \Sigma$, and the generating function for those in the relative interior of $\tau$ is $\prod_{u_i \in \tau} x^{u_i} / (1-x^{u_i}).$  Therefore,
\[
(1-x^{u_1}) \cdots (1-x^{u_s}) \cdot \Hilb(\sigma) =  \sum_{\tau \in \Sigma} \bigg( \prod_{u_i \in \tau} x^{u_i} \prod_{u_j \not \in \tau} (1-x^{u_j}) \bigg).
\]
If $\tau_1, \ldots, \tau_r$ are the maximal cones of $\Sigma$, then taking leading forms at $1_T$ on both sides gives
\[
 u_1 \cdots u_s \cdot (-1)^n \cdot \Hilb(\sigma)_\circ = \sum_{i = 1}^r \,  \prod_{u_j \not \in \tau_i} u_j,
\]
provided that the right hand side is nonvanishing.  Since all of the $u_j$ lie in $\sigma^*$, the right hand side is strictly positive on the interior of $\sigma$.  In particular, it does not vanish, so the degree of $\Hilb(\sigma)_\circ$ is $-n$.
\end{proof}

\noindent Lemma~\ref{degree hilb} can also be seen as a special case of more general results on multigraded Hilbert series of modules.  See \cite[Definition~8.45 and Claim~8.54]{MillerSturmfels05}.

\begin{lemma} \label{unimodular hilb}
Let $\sigma$ be a unimodular cone, spanned by a basis $e_1, \ldots, e_n$ for $N$.  Then the principal part of $\Hilb(\sigma)$ at $1_T$ is
\[
\Hilb(\sigma)_\circ = \frac{(-1)^n}{e_1^* \cdots e_n^*},
\]
where $e_1*, \ldots, e_n^*$ is the dual basis for $M$.
\end{lemma}

\begin{proof}
The generating function $\Hilb(\sigma)$ is given by
\[
\Hilb(\sigma) = \frac{1}{(1-x^{e_1^*}) \cdots (1-x^{e_n^*})}.
\]
Now, $(1-x^{e_i^*})$ is a local parameter at $1_T$, with principal part $(1-x^{e_i^*})_\circ = -e_i^*$.  Since principal parts are multiplicative, it follows that the principal part of $1/(1-x^{e_i^*})$ is $-1/e_i^*$, and the lemma follows.
\end{proof}

\begin{proposition} \label{computing multiplicities}
Let $\sigma$ be an $n$-dimensional rational polyhedral cone in $N_\R$.
\begin{enumerate}
\item If $\sigma_1, \ldots, \sigma_r$ are the maximal cones in a rational polyhedral subdivision of $\sigma$, then
\[
e_\sigma = e_{\sigma_1} + \cdots + e_{\sigma_r}.
\]
\item If $\sigma$ is unimodular, spanned by a basis $e_1, \ldots, e_n$ for $N$, then
\[
e_\sigma = \frac{1}{e_1^* \cdots e_n^*} \ .
\]
\end{enumerate}
\end{proposition}

\noindent In particular, the sum determined by (1) and (2) is independent of the choice of unimodular subdivision.

\begin{proof}
Part (1) follows from the additivity of $\Hilb(\sigma_i)$ (Lemma~\ref{additive hilb}) and the fact that $\Hilb(\sigma)$ and the $\Hilb(\sigma_i)$ all have principal parts in degree $-n$ (Lemma \ref{degree hilb}).  See Proposition \ref{principal degrees}, in the Appendix.  Part (2) is an immediate consequence of Lemma~\ref{unimodular hilb}.
\end{proof}

Recall that for any cone $\tau \in \Delta$,
$\Delta_\tau$ is the fan in $\big(N / (N \cap \Span \tau)\big)_\R$
whose cones are the projections of the cones in $\Delta$ that
contain $\tau$.  If $\sigma$ is a maximal cone containing $\tau$, we define
$e_{\sigma, \tau}$ to be $e_{\overline{\sigma}}$, where
$\overline{\sigma}$ is the image of $\sigma$ in $\Delta_\tau$.  So
$e_{\sigma, \tau}$  is a homogeneous rational function of degree
$(\dim \tau - n)$ in the graded subring $\Sym^\pm(\tau^\perp \cap
M)$ of $\Sym^\pm (M)$.  We write $V(\tau)$ for the $T$-invariant subvariety of $X$ corresponding to $\tau$.

\begin{corollary} \label{compare_to_brion}
If $\sigma$ is an $n$-dimensional rational polyhedral cone in $N_\R$ and $\tau$ is a face of $\sigma$, then
\[
e_\sigma = e_{x_\sigma}[X] \mbox{ \ and \ } e_{\sigma, \tau} = e_{x_\sigma}[V(\tau)].
\]
\end{corollary}

\begin{lemma} \label{relative unimodular}
If $\sigma$ is a unimodular cone spanned by a basis $e_1, \ldots, e_n$ for $N$ and $\tau \preceq \sigma$ then
\[
e_{\sigma, \tau} = \prod_{e_i \not \in \tau} \frac{1}{e_i^*}.
\]
\end{lemma}

\begin{proof} Apply part (2) of Proposition \ref{computing multiplicities} to the fan $\Delta_\tau$.
\end{proof}

\section{Localization and Minkowski weights} \label{localizationsection}

Here we use equivariant multiplicities to describe the natural map
from piecewise polynomials on a complete fan to Minkowski weights.
We then use localization to show that this map agrees with
$\iota^*: A^*_T(X) \rightarrow A^*(X)$.

\begin{lemma} \label{vanishing}
Let $\Delta$ be a complete $n$-dimensional fan.  Then the sum of the rational functions $e_\sigma$ for all maximal cones $\sigma \in \Delta$ is given by
\[
\sum_{\sigma} e_\sigma = \left\{ \begin{array}{ll} 0 & \mbox{ for } n \geq 1. \\ 1 & \mbox{ for } n = 0. \end{array} \right.
\]
\end{lemma}

\begin{proof}
If $n = 0$, then $\Delta$ contains only one cone $0$, and $e_0 = 1$.  Suppose $n \geq 1$. In the polytope algebra ,
\[
\sum_\sigma \, [\sigma] = [N_\R] \pm \mbox{ classes of smaller dimensional cones}.
\]
Applying the linear transformation $\nu^*$ gives
\[
\sum_\sigma \Hilb(\sigma) = 1.
\]
Since each of the principal parts $\Hilb(\sigma)_\circ = (-1)^n \cdot e_\sigma$ is homogeneous of degree $-n$, it follows that the sum of these prinicipal parts must vanish (Appendix, Proposition \ref{principal degrees}), and the lemma follows.
\end{proof}

\begin{lemma} \label{restricted vanishing}
Let $\tau$ be a cone in a complete $n$-dimensional fan $\Delta$.  Then
\[
\sum_{\sigma \succeq \tau} e_{\sigma,\tau} = \left\{ \begin{array}{ll} 0 & \mbox{ for } \dim \tau < n. \\ 1 & \mbox{ for } \dim \tau = n. \end{array} \right.
\]
\end{lemma}

\begin{proof}
Apply Lemma~\ref{vanishing} to the fan $\Delta_\tau$.
\end{proof}

Piecewise polynomials are especially well-behaved on unimodular fans, that is, fans in which each maximal cone is
spanned by a basis for the lattice.  Suppose $\Delta$ is a unimodular fan, and $\rho_1, \ldots, \rho_s$ are the rays of $\Delta$.  Let $v_i$ be the primitive generator of $\rho_i$.  Then there is a unique piecewise linear function $\Psi_i \in \PP^1(\Delta)$ whose values at the primitive generators of the rays are given by the Kronecker delta function
\[
\Psi_i(v_j) = \delta_{ij},
\]
and whose values elsewhere are given by extending linearly on each cone.

Then, for any $k$-dimensional cone $\tau \in \Delta$, we have a piecewise polynomial $\Psi_\tau \in \PP^k(\Delta)$ that vanishes away from $\Star(\tau)$, the union of the cones in $\Delta$ that contain $\tau$, defined by
\[
\Psi_\tau = \prod_{v_i \in \tau} \Psi_i,
\]
and $\PP^*(\Delta)$ is generated by $\{ \Psi_\tau \}_{\tau \in \Delta}$ as a $\Sym^*(M)$-module.

\begin{proof}[Proof of Proposition~\ref{pushforward}]
Since equivariant multiplicities are additive with respect to
subdivisions, we may assume that $\Delta$ is unimodular.  Say
$\rho_1, \ldots, \rho_s$ are the rays of $\Delta$ and $v_i$ is the
primitive generator of $\rho_i$.  Since $\PP^*(\Delta)$ is
generated as a $\Sym^*(M)$-module by the piecewise polynomials
$\Psi_\tau$, it suffices to prove that
\[
\sum e_\sigma \cdot (\Psi_\tau)_\sigma
\]
is in $\Sym^*(M)$ for all $\tau$.  Now, if $\sigma$ is spanned by a basis $e_1, \ldots, e_n$ for $N$ and $\tau \preceq \sigma$, then $(\Psi_\tau)_\sigma = \prod_{v_i \in \tau} e_i^*$.  It then follows from Lemma~\ref{relative unimodular} that
\[
e_\sigma \cdot (\Psi_\tau)_\sigma = \left\{ \begin{array}{ll} e_{\sigma, \tau} & \mbox{ for } \sigma \succeq \tau. \\ 0&  \mbox{ otherwise.} \end{array} \right.
\]
Therefore, by Lemma~\ref{restricted vanishing}, $\sum e_\sigma \cdot (\Psi_\tau)_\sigma$ vanishes unless $\tau$ is a maximal cone, in which case the sum is equal to one.  In particular, $\sum e_\sigma \cdot (\Psi_\tau)_\sigma$ is in $\Sym^*(M)$, as required.
\end{proof}

\begin{proof}[Proof of Proposition~\ref{restricted pushforward}]
The sum in Proposition~\ref{restricted pushforward} is over the maximal cones in $\Star(\tau)$.   Then the proof of Proposition~\ref{restricted pushforward} is similar to the proof of Proposition~\ref{pushforward}, since $\PP^*(\Star(\tau))$ is generated as a $\Sym^*(M)$-module by the restrictions of the piecewise polynomial functions $\Psi_\gamma$, for $\gamma \in \Star(\tau)$.
\end{proof}

It remains to show that if $f$ is a homogeneous piecewise polynomial of degree $k$, then the integer-valued function $c$ on codimension $k$ cones of $\Delta$ given by
\[
c(\tau) = \sum_{\sigma \succeq \tau} e_{\sigma, \tau} f_\sigma
\]
is a Minkowski weight of codimension $k$, and that $f \mapsto c$
agrees with the natural map $\iota^*: A^*_T(X) \rightarrow
A^*(X)$.  Although the entire statement can be proved using the
general machinery of localization, the fact that the integers
$c(\tau)$ give a Minkowski weight is purely combinatorial, and we
include an elementary proof.

We recall the definition of Minkowski weights from
\cite{FultonSturmfels97}.  If $\gamma$ is a codimension $k+1$ cone
in $\Delta$ contained in a codimension $k$ cone $\tau$, we write
$v_{\tau/\gamma} \in N / (N \cap \Span \gamma)$ for the primitive
generator of the image of $\tau$ in $\Delta_\gamma$.

\begin{definition} An integer valued function $c$ on codimension $k$ cones $\tau \in \Delta$ is a Minkowski weight if, for every codimension $k + 1$ cone $\gamma \in \Delta$,
\[
\sum_{\tau \succeq \gamma} c(\tau) \cdot v_{\tau/\gamma} = 0.
\]
\end{definition}

We will use the following basic property of equivariant mulitplicities to show that the integer-valued function $c$ coming from a piecewise polynomial function is a Minkowski weight.  Let $v_\rho$ denote the primitive generator of a ray $\rho$.

\begin{proposition} \label{divisor relation}
If $\sigma$ is an $n$-dimensional rational polyhedral cone in $N_\R$ and $u$ is in $M$ then
\[
\sum_{\rho \preceq \sigma} \<u, v_\rho \> \cdot e_{\sigma,\rho} = u \cdot e_\sigma.
\]
\end{proposition}

\noindent We will prove the proposition by subdividing $\sigma$ and reducing to the case where $\sigma$ is unimodular.

\begin{lemma} \label{unimodular prop}
If $\sigma$ is an $n$-dimensional unimodular cone in $N_\R$ and $u$ is in $M$ then
\[
\sum_{\rho \preceq \sigma} \<u, v_\rho\> \cdot e_{\sigma, \rho} = u \cdot e_\sigma.
\]
\end{lemma}

\begin{proof}
Say $\sigma$ is spanned by a basis $e_1, \ldots, e_n$ for $N$, and $u = u_1 e_1^* + \cdots u_n e_n^*$.  Then
\[
\sum \<u,v_\rho\> \cdot e_{\sigma, \rho} = \sum_{i=1}^n  \frac{u_i}{e_1^* \cdots \widehat{e_i^*} \cdots e_n^*},
\]
which is equal to $u \cdot e_\sigma$.
\end{proof}

\begin{lemma} \label{ray sum}
If $\sigma_1, \ldots, \sigma_s$ are the maximal cones in a subdivision of $\sigma$, and if $\rho$ is a ray in this subdivision then
\[
\sum_{\sigma_i \succeq \rho} e_{\sigma_i, \rho} =  \left\{ \begin{array}{ll} e_{\sigma, \rho} & \mbox{ if } \rho \preceq \sigma. \\
0 & \mbox{ otherwise.} \end{array} \right.
\]
\end{lemma}

\begin{proof}
Suppose $\rho$ lies in the relative interior of a $k$-dimensional face $\tau \preceq \sigma$.  Consider the fan $\Delta_{\rho}$, whose maximal cones are the images of the $\sigma_i \succeq \rho$.  The support $|\Delta_{\rho}|$ is a closed polyhedral cone in an $(n-1)$-dimensional vector space whose minimal face is $k-1$-dimensional, so the polar dual $|\Delta_{\rho}|^*$ has dimension $n - k$.  It follows that the principal part of $\nu^*(|\Delta_\rho|)$ has degree $k -n$.  Since each $e_{\sigma,\rho}$ has degree $1-n$, and $\sum e_{\sigma_i, \rho}$ is the principal part of $\pm \nu^*(|\Delta_\rho|)$ unless this sum vanishes (Appendix, Proposition \ref{principal degrees}), the lemma follows.
\end{proof}

\begin{proof}[Proof of Proposition~\ref{divisor relation}]
Let $\sigma_1, \ldots, \sigma_r$ be the maximal cones of a unimodular subdivision of $\sigma$.  Then $u \cdot e_\sigma = u \cdot e_{\sigma_1} + \cdots + u \cdot e_{\sigma_r}$.  Since $\sigma_i$ is unimodular,
\[
u \cdot e_{\sigma_i} = \sum_{\rho \preceq \sigma_i} \< u, v_\rho \> e_{\sigma_i, \rho}.
\]
Therefore, by rearranging terms in the summation, we have
\[
u \cdot e_\sigma = \sum_\rho \bigg( \sum_{\sigma_i \succeq \rho} \< u, v_\rho \> e_{\sigma_i, \rho}\bigg).
\]
By Lemma~\ref{ray sum}, the right hand side is equal to $\sum_{\rho \preceq \sigma} e_{\sigma, \rho}$, as required.
\end{proof}

\begin{proposition}
Let $f \in \PP^k(\Delta)$ be a homogeneous piecewise polynomial of degree $k$.  Then the integers
\[
c(\tau) = \sum_{\sigma \succeq \tau} e_{\sigma,\tau} f_\sigma
\]
are a Minkowski weight of codimension $k$ on $\Delta$.
\end{proposition}

\begin{proof}
Let $\gamma$ be a codimension $k + 1$ cone in $\Delta$.  It will suffice to show that $\sum \< u, v_{\tau/\gamma} \> c(\tau) = 0$ for any $u \in (M \cap \gamma^\perp)$, where the sum is over all codimension $k$ cones $\tau$ containing $\gamma$.  To prove this, we will show that $\sum \<u, v_{\tau/\gamma} \> c(\tau)$, which is an integer by Proposition~\ref{restricted pushforward}, is divisible by the linear function $u$ in $\Sym^*(M)$.  Now,
\[
\sum_{\tau} \<u,v_{\tau/\gamma} \> c(\tau) = \sum_{\tau} \bigg( \sum_{\sigma \succeq \tau} \<u, v_{\tau/\gamma} \> e_{\sigma, \tau} f_\sigma\bigg),
\]
and the sum on the right hand side may be rearranged as
\[
\sum_\sigma \bigg( f_\sigma \cdot \sum_{\tau \preceq \sigma} \< u, v_{\tau/\gamma} \> e_{\sigma, \tau} \bigg).
\]
Applying Proposition~\ref{divisor relation} to $\Delta_\gamma$ then gives
\[
\sum_{\tau \preceq \sigma} \< u, v_{\tau/\gamma} \> e_{\sigma, \tau} = u \cdot e_{\sigma,  \gamma}.
\]
It follows that the integer $\sum \<u, v_{\tau/\gamma} \> c(\tau)$ is divisible by $u$ in $\Sym^*(M)$, as claimed, and hence must vanish.
\end{proof}

\begin{proof}[Proof of Theorem \ref{localization}]
To show that $f \mapsto c$ agrees with the natural map $\iota^*:
A^*_T(X) \rightarrow A^*(X)$, we must prove that
\[
\int_{V(\tau)} \iota^*c_f  = c(\tau),
\]
where $c_f$ denotes the equivariant Chow cohomology class whose
restriction to a torus fixed point is $f_\sigma \in \Sym^*(M)
\cong A^*_T(x_\sigma)$.  By Corollary \ref{compare_to_brion},
$e_{\sigma, \tau}$ is equal to the equivariant multiplicity
$e_{x_\sigma}[V(\tau)]$ of the nondegenerate $T$-fixed point
$x_\sigma$ in $V(\tau)$.  Therefore, by localization
\cite{EdidinGraham98b}, in equivariant Chow homology tensored with
$\Sym^\pm(M)$, we have
\[
\int_{V(\tau)} c_f  = \sum_\sigma e_\sigma f_\sigma.
\]
Since $\sum_\sigma e_\sigma f_\sigma$ is an integer,
projecting to $A_*(X)$ gives $\int_{V(\tau)} \iota^* c_f  = c(\tau)$.
\end{proof}

\section{Applications to Chow cohomology of toric varieties} \label{applications}

Here we use combinatorial computations with piecewise polynomials to study the map $\iota^* : A^*_T(X) \rightarrow A^*(X)$ for some specific complete toric varieties $X$.  As discussed in the introduction, this map is known to be surjective with kernel generated by $M$ in degree one if $X$ is smooth, and similar statements hold over $\Q$ if $X$ is simplicial.  We give the first examples showing that $\iota^*$ is not surjective in general, and that its kernel is not always generated in degree one.

\begin{example}[Mirror dual of $\P^1 \times \P^1$] \label{mod Z2}
Let $N = \Z^2$, and let $\Delta$ be the complete fan in $\R^2$ whose rays are generated by
\[
\begin{array}{llll}
v_1 = (1,1),   &  v_2 = (1,-1),  &  v_3 = (-1,-1),  &  v_4 = (-1,1),
\end{array}
\]
and whose maximal cones are
\[
\begin{array}{llll}
\sigma_1 = \< v_1, v_2 \>, & \sigma_2 = \<v_2, v_3 \>, & \sigma_3 = \< v_3, v_4 \>, & \sigma_4 = \< v_1, v_4 \>.
\end{array}
\]
Then $X = X(\Delta)$ is isomorphic to $(\P^1 \times \P^1) / \Z_2$, which is the Fano surface that is ``mirror dual" to $\P^1 \times \P^1$.

We claim that the image of $\PP^2(X)$ under the map $f \mapsto
\sum_{i=1}^4 e_{\sigma_i} f_{\sigma_i}$ is exactly $2 \Z$.  We
compute $e(\sigma_1)$ using the unimodular subdivision of
$\sigma_1$ along $v = (1,0)$,
\[
\sigma_1 = \< v_1, v \> \cup \< v, v_2 \>.
\]
Then, writing $a = e_1^*$ and $b = e_2^*$, the dual cones of $\<v_1, v\>$ and $\< v, v_2 \>$ are $\< b, a-b \>$ and $\<b, a+b\>$, respectively, so
\begin{eqnarray*}
e(x_1) & = & \frac{1}{b (a-b)} - \frac{1}{b (a+b)}. \\
  & = & \frac{2}{a^2 - b^2}.
\end{eqnarray*}

\noindent Similarly, we compute $e(\sigma_3) =2 / (a^2-b^2)$ and
\[
e(\sigma_2) = e(\sigma_4) = \displaystyle{\frac{-2}{a^2 - b^2}}.
\]

Therefore, since two divides every term in $\sum_{i = 1}^4 e_{\sigma_i} f_{\sigma_i}$, the sum must be divisible by two.  Also, the piecewise polynomial function $f$ that vanishes on $\sigma_2 \cup \sigma_3 \cup \sigma_4$ and whose restriction to $\sigma_1$ is $a^2 - b^2$ maps to two.  So the image of $\PP^2(\Delta)$ is $2 \Z$, as required.
\end{example}

\begin{proof}[Proof of Theorem~\ref{surface}]
Applying Theorem~\ref{localization} to Example~\ref{mod Z2} shows that the image of $\iota^*: A^2_T(X) \rightarrow A^2(X)$ is $2A^2(X)$, which is a proper subgroup of $A^2(X) \cong \Z$.
\end{proof}

In the following examples, we consider fans in $\R^3$ with respect to the lattice $N = \Z^3$.

\begin{example}[Mirror dual of $\P^1 \times \P^1 \times \P^1$] \label{first}
Consider the toric variety $X = X(\Delta)$, where $\Delta$ is the fan whose nonzero cones are the cones over the faces of the cube with vertices $(\pm1, \pm1, \pm1)$.  Then $X$ is the Fano toric threefold that is ``mirror dual" to $\P^1 \times \P^1 \times \P^1$.  Recall that, since $X$ is complete, the rank of $A^i(X)$ is equal to the rank of $A_i(X)$ \cite[Proposition~2.4]{FultonSturmfels97}, so $\rk A^0(X) = \rk A^3(X) = 1$.  Furthermore, since $A_2(X)$ is the Weil divisor class group of $X$, we also have $\rk A^2(X) = 5$.  The remainder of the following table can be filled in by straightforward linear algebra computations with piecewise polynomial functions.

\vspace{10 pt}
\begin{center}

\begin{tabular}{|r|c|c|c|} \hline
& $ \rk A^i_T(X)$ & $\rk M \cdot A^{i-1}_T(X)$ & $\rk A^i(X)$ \\ \hline
$i=0$ & $1$ &$ 0$ & $1$ \\
$i=1$ & $4$ & $3$ & $1$ \\
$i=2$ & $11$ & $9$ & $5$ \\
$i=3$ & $23$ & $22$ & $1$ \\ \hline
\end{tabular}

\vspace{10 pt}
\end{center}

\noindent From these computations, it is clear that $A^2_T(X)$ does not surject onto $A^2(X)$, since its image has rank at most two.
\end{example}

\begin{example}[Fulton's threefold] \label{Fulton}
Consider the toric variety $X' = X(\Delta')$, where $\Delta'$ is the fan combinatorially equivalent to the fan over the cube as in the preceding example, but with the ray through $(1,1,1)$ replaced by the ray through $(1,2,3)$.  Then $X'$ is complete and, as in the previous example, $\rk A^0 (X') = \rk A^3(X') = 1$, and $\rk A^2(X') = 5$, but Fulton showed that $X'$ has no nontrivial line bundles \cite[pp.\ 25--26]{Fulton93}, so $A^1(X') = 0$.  The remainder of the following table is filled in by linear algebra computations with piecewise polynomial functions.

\vspace{10 pt}
\begin{center}

\begin{tabular}{|r|c|c|c|} \hline
& $ \rk A^i_T(X')$ & $\rk M \cdot A^{i-1}_T(X')$ & $\rk A^i(X')$ \\ \hline
$i=0$ & $1$ &$ 0$ & $1$ \\
$i=1$ & $3$ & $3$ & $0$ \\
$i=2$ & $8$ & $6$ & $5$ \\
$i=3$ & $20$ & $16$ & $1$ \\ \hline
\end{tabular}

\vspace{10 pt}
\end{center}

\noindent Here, again, we see that $\iota^*: A^2_T(X') \rightarrow A^2(X')$ is not surjective, since its image has rank at most two.  Furthermore, the kernel of $\iota^*$ is not generated in degree one, since the degree one part of the kernel is $M$, and $A^3_T(X') / M \cdot A^2_T(X')$ has rank four, and hence cannot map injectively into $A^3(X')$.  However, $X'$ is not projective, so to prove Theorem~\ref{threefold}, it remains to give a projective example with similar properties.
\end{example}

\begin{example} \label{final threefold}
Consider the toric variety $X'' = X(\Delta'')$, where $\Delta''$ is the fan combinatorially equivalent to the fan over the cube as in Example \ref{first}, but with the ray through $(1,1,1)$ replaced by the ray through $(1,1,2)$ and with the ray through $(1,-1,1)$ replaced by the ray through $(1,-1,2)$.  It is straightforward to check that $-3K_{X''}$ is Cartier and ample, so $X''$ is $\Q$-Fano and projective.  We compute the following table as in the preceding examples.

\vspace{10 pt}
\begin{center}

\begin{tabular}{|r|c|c|c|} \hline
& $ \rk A^i_T(X'')$ & $\rk M \cdot A^{i-1}_T(X'')$ & $\rk A^i(X'')$ \\ \hline
$i=0$ & $1$ &$ 0$ & $1$ \\
$i=1$ & $4$ & $3$ & $1$ \\
$i=2$ & $10$ & $9$ & $5$ \\
$i=3$ & $22$ & $19$ & $1$ \\ \hline
\end{tabular}

\vspace{10 pt}
\end{center}
\end{example}

\begin{proof}[Proof of Theorem~\ref{threefold}]
From the computations in Example~\ref{final threefold}, we conclude that $\iota^*:A^*_T(X)_\Q \rightarrow A^*(X)_\Q$ is not surjective in degree two, and its kernel in degree three is not in the ideal generated by its kernel in degree one.
\end{proof}

To balance these negative results, we conclude by proving a
positive statement: $\iota^* : A^*_T(X) \rightarrow A^*(X)$ is
always surjective in degree one.

\begin{theorem}  \label{pic}
For any toric variety $X= X(\Delta)$, $\iota^*:A^1_T(X) \rightarrow A^1(X)$ is surjective, giving a natural isomorphism $A^1(X) \cong \PP^1(\Delta) / M$.
\end{theorem}

\begin{proof}
If $X$ is smooth, then the statement is clear.  Suppose $X$ is singular, and let
\[
X_r \rightarrow \cdots \rightarrow X_1 \xrightarrow{\pi} X_0 = X
\]
be a resolution of singularities, where each $X_i = X(\Delta_i)$ is a toric variety and $X_{i+1} \rightarrow X_i$ is the blowup along a smooth $T$-invariant center.  Say $X_1$ is the blowup of $X$ along $V(\tau)$ and $V(\rho) \subset X_1$ is the exceptional divisor.  By induction on $r$, we may assume $A^1(X_1) \cong \PP^1(\Delta_1) / M$.  Also, we may assume $A^1(V(\rho)) = \PP^1(\Star(\rho)) / M$ and $A^1(V(\tau)) = \PP^1(\Star(\tau))/M$, by induction on dimension.  Then $\pi^*: A^1(X) \rightarrow A^1(X_1)$ is injective, and $c \in A^1(X_1)$ is in the image of $\pi^*$ if and only if $c|_{V(\rho)}$ is in the image of $A^1(V(\tau))$ \cite[Theorem~3.1]{Kimura92}.  The theorem then follows, since $\Star(\rho)$ is a subdivision of $\Star(\tau)$, $\Delta_1$ and $\Delta$ coincide everywhere else, and the class of a piecewise linear function $[\Psi] \in \PP^1(\Star(\rho))/M$ is pulled back from $\Star(\tau)$ if and only if $\Psi$ is given by a single linear function on each cone of $\Star(\tau)$.
\end{proof}

\begin{corollary} \label{picard}
For any toric variety $X$, the canonical map $\Pic(X) \rightarrow A^1(X)$ is an isomorphism.
\end{corollary}

\begin{proof}
The corollary follows from the canonical identification of $\PP^1(X) / M$ with $\Pic(X)$ \cite[pp.\ 65--66]{Fulton93}.
\end{proof}

\noindent Corollary \ref{picard} was known previously in the case where $X$ is complete \cite{Brion89}.  See also \cite[Corollary~3.4]{FultonSturmfels97}.

\begin{remark}
One can use Kimura's inductive method, as in the proof of Theorem~\ref{picard} and \cite[Theorem~1]{chowcohom}, to compute the Chow cohomology of an arbitrary toric variety in all degrees.  However, the resulting induction is more subtle, as Theorems~\ref{surface} and \ref{threefold} suggest.
\end{remark}

\section{Localization formula for mixed volumes of lattice polytopes} \label{mixed volumes}

Let $P_1, \ldots, P_n$ be lattice polytopes in $M_\R$.  For nonnegative real numbers $a_i$, the euclidean volume of $a_1P_1 + \cdots + a_n P_n$ is a homogeneous polynomial function of $(a_1, \ldots, a_n)$.  The \emph{mixed volume} $V(P_1, \ldots, P_n)$ is defined to be the coefficient of $a_1 \cdots a_n$ in this polynomial.  Let $\Delta$ be the inner normal fan to $P_1 + \cdots +P_n$, and let $u_i(\sigma) \in M$ be the vertex of $P_i$ that is minimal on $\sigma$, for each maximal cone $\sigma \in \Delta$.

\begin{theorem} \label{mixed volumes formula}
The mixed volume of the polytopes $P_i$ is given by
\[
n! \cdot V(P_1, \ldots, P_n) = (-1)^n \sum_{\sigma \in \Delta} e_\sigma \cdot u_1(\sigma) \cdots u_n(\sigma).
\]
\end{theorem}

\noindent Theorem~\ref{mixed volumes formula} follows from Theorem~\ref{localization} and the fact that $V(P_1, \ldots, P_n)$ is the degree of $D_1 \cdots D_n$, where $D_i$ is the $T$-Cartier divisor on $X(\Delta)$ corresponding to $P_i$ \cite[p.~116]{Fulton93}.  However, the statement of the theorem is purely combinatorial, and we give a combinatorial proof based on Brion's formula for generating functions for lattice points in polyhedra.  The methods used in this proof may be of independent interest.

Let $P$ be a lattice polytope in $M_\R$.  Let $\Delta$ be the normal fan to $P$, and let $u(\sigma) \in M$ be the vertex of $P$  that is minimal on $\sigma$, for each maximal cone $\sigma \in \Delta$.

\begin{BF}
The generating function for lattice points in $P$ is
\[
\sum_{u \in (P \cap M)} x^u = \sum_\sigma x^{u(\sigma)} \cdot \Hilb(\sigma).
\]
\end{BF}

In addition to Brion's Formula, we will use the following formula for mixed volumes, which is a lattice point counting analogue of the alternating sum of volumes in formula (3) of \cite[p.\ 116 ]{Fulton93}.

\begin{proposition} \label{lattice point formula}
Let $P_1, \ldots, P_n$ be lattice polytopes in $M_\R$.  Then
\[
n! \cdot V(P_1, \ldots, P_n) = \sum_{1 \leq i_1 < \cdots < i_k \leq n} (-1)^{n-k} \#\big( (P_{i_1} + \cdots + P_{i_k}) \cap M \big).
\]
\end{proposition}

\begin{proof}
The number of lattice points in $a_1P_1 + \cdots + a_n P_n$ is a polynomial in the $a_i$ of degree at most $n$, and the degree $n$ part of this polynomial is $n!$ times the volume of $a_1 P_1 + \cdots + a_n P_n$  \cite[Theorem~7]{McMullen78}.  Therefore, $n! \cdot V(P_1, \ldots, P_n)$ is the coefficient of $a_1 \cdots a_n$ in this polynomial, and the proposition is an immediate consequence of the following lemma.
\end{proof}

\begin{lemma}
Let $f \in \R[t_1, \ldots, t_n]$ be a polynomial function on $\R^n$ of degree at most $n$.  The coefficient of $t_1 \cdots t_n$ in $f$ is
\[
\sum_{1 \leq i_1 < \cdots < i_k \leq n} (-1)^{n-k} f(e_{i_1} + \cdots + e_{i_k})
\]
where $\{e_1, \ldots, e_n\}$ is the standard basis for $\R^n$.
\end{lemma}

\begin{proof}
The function taking a polynomial $g$ to $\sum (-1)^{n-k} g(e_{i_1} + \cdots + e_{i_k})$ vanishes on any monomial that does not contain all $n$ variables, and its value on $t_1 \cdots t_n$ is 1.
\end{proof}

\begin{proof}[Proof of Theorem \ref{mixed volumes formula}]
For each $\sigma$, $(-1)^n \cdot e_\sigma \cdot u_1(\sigma) \cdots u_n(\sigma)$ is the principal part of
\[
(x^{u_1(\sigma)}-1) \cdots (x^{u_n(\sigma)} - 1) \cdot \Hilb(\sigma).
\]
Expanding the product of the binomials, taking the sum over all $\sigma$, and applying Brion's Formula then gives
\[
\sum_{1 \leq i_1 < \cdots < i_k \leq n} (-1)^{n-k} \cdot \sum_{u \in (P_{i_1} + \cdots + P_{i_k} \cap M)} x^u.
\]
The theorem then follows from Proposition~\ref{lattice point formula} by taking principal parts, since the leading form of $x^u$ at $1_T$ is equal to one.
\end{proof}

\section{Bott residue formula for toric vector bundles}

The mixed volume $V(P_1, \ldots, P_n)$ is the degree of the top Chern class of the toric vector bundle $\O(D_1) \oplus \cdots \oplus \O(D_n)$, where $D_i$ is the $T$-Cartier divisor corresponding to $P_i$.  Therefore, mixed volumes are a special case of Chern numbers of toric vector bundles, and Theorem~\ref{mixed volumes formula} may be generalized as follows.  Given a multiset of linear functions $\bu \subset M$ let $\varepsilon_i(\bu) \in \Sym^i(M)$ be the $i$-th elementary symmetric function in the elements of $\bu$.  For instance, if $\bu = \{ u_1, \ldots, u_r \}$, then $\varepsilon_1(\bu) = u_1 + \cdots + u_r$ and $\varepsilon_r (\bu) = u_1 \cdots u_r$.  For a partition $\lambda = (\lambda_1, \ldots, \lambda_s)$ of $n$, let $\varepsilon_\lambda(\bu) \in \Sym^n(M)$ be the product
\[
\varepsilon_\lambda(\bu) = \varepsilon_{\lambda_1}(\bu) \cdots \varepsilon_{\lambda_s} (\bu).
\]

Recall that, for any toric vector bundle $\E$ on an arbitrary toric variety $X = X(\Delta)$ and any maximal cone $\sigma \in \Delta$, there is a unique multiset $\bu(\sigma) \subset M$ such that the restriction of $\E$ to $U_\sigma$ splits equivariantly as
\[
\E|_{U_\sigma} \cong \bigoplus_{u \in \bu(\sigma)} \O(\divisor \chi^{u}).
\]
See \cite{Klyachko90} or \cite[Section~2]{branchedcovers} for this and other basic facts about toric vector bundles.

\begin{theorem} \label{chern numbers}
Let $\E$ be a toric vector bundle on a complete toric variety $X$, and let $\lambda$ be a partition of $n$.  Then the Chern number $c_\lambda(\E)$ is given by
\[
c_\lambda(\E) = \sum_\sigma e_\sigma \cdot \varepsilon_\lambda(\bu(\sigma)).
\]
\end{theorem}

\begin{proof}
The Chern number $c_\lambda(\E)$ is equal to the integral over $[X]$ of the equivariant Chow cohomology class corresponding to the piecewise polynomial whose restriction to $\sigma$ is $\varepsilon_\lambda(\bu(\sigma))$.  Therefore, the theorem follows from Theorem~\ref{localization}.
\end{proof}

Theorem~\ref{chern numbers} has a straightforward generalization to top degree polynomials in the Chern classes of several toric vector bundles (we omit the details), and may be seen as a Bott residue formula for vector bundles on toric varieties with arbitrary singularities.  This solves the toric case of the problem of proving residue formulas on singular varieties posed in \cite[Section~5]{EdidinGraham98b}.  Edidin and Graham handled the case of toric subvarieties of smooth toric varieties, but the extension to arbitrary toric varieties is nontrivial; there are singular toric varieties, such as Fulton's threefold (Example~\ref{Fulton}) that have no nontrivial line bundles, and hence admit no nonconstant morphisms to smooth varieties.

\begin{example}
We apply Theorem~\ref{chern numbers} to compute the Chern numbers of a specific nonsplit rank two toric vector bundle on the singular toric variety $X = X(\Delta)$ mirror dual to $\P^1 \times \P^1 \times \P^1$ (see Example~\ref{first}).  For this example, we assume that the base field has at least three elements. The primitive generators of the rays of $\Delta$ are
\[
\begin{array}{llll}
v_1 = ( 1,1,1 ) , &  v_2 = (1,1,-1), & v_3 = (1,-1,1), & v_4 = ( 1,-1,-1 ) \\
v_5 = (-1,1,1), & v_6 = ( -1,1,-1 ), & v_7 = ( -1,-1,1 ), & v_8 = (-1,-1,-1),
\end{array}
\]
and the maximal cones of $\Delta$ are
\[
\begin{array}{lll}
\sigma_1 = \<v_1,v_2,v_3,v_4\>, & \sigma_2 = \<v_1,v_2,v_5,v_6\>, & \sigma_3 = \<v_1,v_3,v_5,v_7\>, \\
\sigma_4 = \<v_2,v_4,v_7,v_8\>, & \sigma_5 = \<v_3,v_4,v_7,v_8\>, & \sigma_6 = \<v_5,v_6,v_7,v_8\>.
\end{array}
\]

Let $\rho_i$ be the ray of $\Delta$ spanned by $v_i$, let $E = k^2$, fix four distinct lines $L_1$, $L_2$, $L_3$, and $L_4$ in $E$, and let $\E$ be the toric vector bundle determined by the filtrations
\[
\begin{array}{ll}
\begin{array}{ll}
{ E^{\rho_1}(i) = \left\{ \begin{array}{ll} E & \mbox{ for } i  \leq -1,  \\ L_1 & \mbox{ for } 0 \leq i \leq 3, \\ 0 & \mbox{ for } i  > 3, \end{array} \right. } &
{ E^{\rho_4}(i) = \left\{ \begin{array}{ll} E & \mbox{ for } i  \leq -1,  \\ L_2 & \mbox{ for } 0  \leq i   \leq 3, \\ 0 & \mbox{ for } i  > 3, \end{array} \right. } \\ \\
{ E^{\rho_6}(i) = \left\{ \begin{array}{ll} E & \mbox{ for } i  \leq -1,  \\ L_3 & \mbox{ for } 0 \leq i \leq 3, \\ 0 & \mbox{ for } i  > 3, \end{array} \right. } &
{ E^{\rho_7}(i) = \left\{ \begin{array}{ll} E & \mbox{ for } i  \leq -1,  \\ L_4 & \mbox{ for } 0  \leq i   \leq 3, \\ 0 & \mbox{ for } i  > 3, \end{array} \right. } \\ \\
\end{array}
\end{array}
\]
and
\[
E^{\rho_j}(i) = \left\{ \begin{array}{ll} E & \mbox{ for } i  \leq 1, \\ 0 & \mbox{ for } i  > 1, \end{array} \right.
\]
for $j \in \{ 2,3,5, 8 \}$.  Since the lines $L_i$ are distinct, the vector bundle $\E$ does not split as a sum of line bundles.  It is straightforward to check that the multisets of linear functions $\bu(\sigma_i)$ determined by $\E$ are as follows.  For simplicity, we write $a$, $b$, and $c$, for $e_1^*$, $e_2^*$ and $e_3^*$, respectively.
\[
\begin{array}{ll}
\bu(\sigma_1) = \{ (a+b+c), (a-b-c) \}, & \bu(\sigma_2) = \{(a+b+c),(-a+b-c)\}, \\
\bu(\sigma_3) = \{(a+b+c),(-a-b+c)\}, & \bu(\sigma_4) = \{(a-b-c),(-a+b-c)\}, \\
\bu(\sigma_5) = \{(a-b-c),(-a-b+c)\}, & \bu(\sigma_6) = \{(-a+b-c),(-a-b+c)\}.
\end{array}
\]

To compute the Chern numbers of $\E$, we now need only to compute the equivariant multiplicities $e_{\sigma_i}$.  First, $\sigma_1^*$ is spanned by $u_1= (1,1,0)$, $u_2 = (1,0,1)$, $u_3 = (1,0,-1)$, and $u_4 = (1,-1,0)$.  Let $u = (1,0,0)$.  We compute, as in \cite[Example~1.8]{ehrhartseries},
\[
\Hilb(\sigma_1) = \frac{(1+x^u)(1-x^{2u})}{(1-x^{u_1})(1-x^{u_2})(1-x^{u_3})(1-x^{u_4})}.
\]
Since the principal parts of $1+x^u$, $1-x^{2u}$ and $1-x^{u_i}$ at $1_T$ are $2$, $-2u$, and $-u_i$ respectively, it is then straightforward to compute $e_{\sigma_1} = - \Hilb(\sigma)_\circ$.  Then
\[
e_{\sigma_1} = \frac{4a}{(b^2-a^2)(c^2-a^2)}.
\]
By symmetry, $e_{\sigma_6} = -e_{\sigma_1}$, and similarly
\[
e_{\sigma_2}  = \frac{4b}{(a^2-b^2)(c^2-b^2)} = -e_{\sigma_5},
\]
and
\[
e_{\sigma_3} = \frac{4c}{(a^2-c^2)(b^2-c^2)} = -e_{\sigma_4} .
\]
Then, using Theorem~\ref{chern numbers} and combining the summands coming from $\sigma_i$ and $\sigma_{7-i}$,
\vspace{5 pt}
\[
c_{111}(\E) =  \frac{2 \cdot (2a)^3 \cdot 4a}{(b^2-a^2)(c^2-a^2)} + \frac{2 \cdot (2b)^3 \cdot 4b}{(a^2-b^2)(c^2-b^2)} + \frac{2 \cdot (2c)^3 \cdot 4c}{(a^2-c^2)(b^2-c^2)},
\]

\vspace{5 pt}
\noindent which simplifies to $c_{111}(\E) = 64$.  Similarly,
\vspace{5 pt}
\[
c_{21}(\E) = \frac{16a^2 (a^2 - b^2 - c^2)}{(b^2-a^2)(c^2-a^2)} + \frac{16b^2(-a^2 + b^2 -c^2)}{(a^2-b^2)(c^2-b^2)} + \frac{16c^2(-a^2 - b^2 + c^2)}{(a^2-c^2)(b^2-c^2)},
\]

\vspace{5 pt}
\noindent which simplifies to $c_{21}(\E) = 32$.
\end{example}

\section{Appendix: Principal parts of rational functions}

Associated graded rings and leading forms have been standard tools for about as long as commutative algebra has been applied to local algebraic geometry \cite{Samuel53, Samuel55}.  The generalization from leading forms of regular functions to principal parts of rational functions is straightforward but, since we have been unable to locate a reference, we include a brief account.

Let $X$ be an algebraic variety over a field $k$, and let $x \in X(k)$ be a smooth point.  Let $\m$ be the maximal ideal in the local ring $\O_{X,x}$.  Since $x$ is smooth,
\[
(\m^d / \m^{d+1}) \cong \Sym^d(\m/\m^2),
\]
for all nonnegative integers $d$ \cite[Theorem 11.22]{AtiyahMacdonald69}.  Suppose $g \in \O_{X,x}$ is a regular function whose order of vanishing at $x$ is $d$.  Then the \emph{leading form} of $g$ is its image
\[
g_\circ \in \Sym^d(\m / \m^2).
\]
In other words, if $x_1, \ldots, x_n$ is a local system of parameters, then $g$ can be expanded uniquely as a power series in $k[[x_1, \ldots, x_n]]$, and the sum of the lowest degree terms in this power series is the homogeneous degree $d$ polynomial in $x_1, \ldots, x_n$ that maps to $g_\circ$ under the canonical isomorphism
\[
k[x_1, \ldots, x_n]_d \cong \Sym^d (\m / \m^2).
\]

Now $\m / \m^2$ is the cotangent space of $X$ at $x$, so $g_\circ$ is naturally a regular function on the tangent space $T_{X,x}$, and the zero locus of $g_\circ$ is the tangent cone of the divisor of zeros of $g$ at $x$ \cite[Lecture~20]{Harris92}.  Note that leading forms are multiplicative; if $g,h \in \O_{X,x}$, then
\[
(gh)_\circ = g_\circ h_\circ
\]
in $\Sym^*(\m/\m^2)$.  For convenience, we define the leading form of zero to be $0 \in \Sym^*(\m/\m^2)$.

\vspace{ 5 pt}

Suppose $f$ is a rational function on $X$.  Then $f$ can be written as a fraction $f = g/h$, with $g, h \in \O_{X,x}$.  We define the \emph{principal part} of $f$ to be
\[
f_\circ = g_\circ / h_\circ,
\]
which is a homogeneous element of $\Sym^\pm(\m/\m^2)$, the $\Z$-graded ring obtained by inverting all homogeneous elements in $\Sym^*(\m/\m^2)$.  Note that $f_\circ$ is well-defined; if $g/h = g'/h'$, then $gh' = g'h$ (since $\O_{X,x}$ is a domain), and therefore
\[
g_\circ h'_\circ = g'_\circ h_\circ,
\]
since leading forms are multiplicative, so $g_\circ/h_\circ = g'_\circ / h'_\circ$.  Also, $f_\circ$ is naturally a rational function on $T_{X,x}$, and its divisors of zeros and poles are the tangent cones of the zeros and poles of $f$, respectively.

\begin{proposition} \label{principal degrees}
Suppose $f_1, \ldots, f_s$ are rational functions on $X$ with principal parts in degree $d$, and let $f = f_1 + \cdots + f_s$.  Then either
\[
f_\circ = (f_1)_\circ + \cdots + (f_s)_\circ,
\]
or $(f_1)_\circ + \cdots + (f_s)_\circ = 0$ and the principal part of $f$ is in degree strictly greater than $d$.
\end{proposition}

\begin{proof}
If $f = 0$ then the proposition is clear.  Suppose $f$ is nonzero, and express each $f_i$ as a fraction $f_i = g_i/h_i$, with $g_i,h_i \in \O_{X,x}$.  Then we can write $f$ as a fraction over a common denominator
\[
f \ = \ \sum_{i=1}^s \frac{g_i \cdot h_1 \cdots \widehat{h_i} \cdots h_s}{h_1 \cdots h_s} \, .
\]
Say $h_i$ vanishes to order $d_i$ at $x$.  Then each summand in the numerator above vanishes to order exactly $d_1 + \cdots + d_s  + d$.  Therefore, either the numerator vanishes to order exactly $d_1 + \cdots + d_s + d$ and $f_\circ = (f_1)_\circ + \cdots + (f_s)_\circ$, or the numerator vanishes to some larger order and $f_\circ$ has degree greater than $d$.
\end{proof}

\begin{corollary} \label{regular principal parts}
Suppose $f_1, \ldots, f_s$ are rational functions on $X$ with principal parts in degree $d$, and suppose $f = f_1 + \cdots + f_s$ is regular at $x$.  Then $(f_1)_\circ + \cdots + (f_s)_\circ \in \Sym^*(\m / \m^2)$ is regular on $T_{X,x}$.
\end{corollary}

\begin{proof}
By Proposition~\ref{principal degrees}, if $(f_1)_\circ + \cdots + (f_s)_\circ$ does not vanish then it is equal to $f_\circ$, which is the principal part of a regular function.
\end{proof}

\bibliography{math}

\providecommand{\bysame}{\leavevmode\hbox to3em{\hrulefill}\thinspace}
\providecommand{\MR}{\relax\ifhmode\unskip\space\fi MR }
\providecommand{\MRhref}[2]{%
  \href{http://www.ams.org/mathscinet-getitem?mr=#1}{#2}
}
\providecommand{\href}[2]{#2}
\begin{thebibliography}{McM79}

\bibitem[AM69]{AtiyahMacdonald69}
M.~Atiyah and I.~Macdonald, \emph{Introduction to commutative algebra},
  Addison-Wesley Publishing Co., Reading, Mass.-London-Don Mills, Ont., 1969.

\bibitem[Bar02]{Barvinok02}
A.~Barvinok, \emph{A course in convexity}, Graduate Studies in Mathematics,
  vol.~54, American Mathematical Society, Providence, RI, 2002.

\bibitem[Bri89]{Brion89}
M.~Brion, \emph{Groupe de {P}icard et nombres caract\'eristiques des
  vari\'et\'es sph\'eriques}, Duke Math. J. \textbf{58} (1989), no.~2,
  397--424.

\bibitem[Bri96]{Brion96}
\bysame, \emph{Piecewise polynomial functions, convex polytopes and enumerative
  geometry}, Parameter spaces (Warsaw, 1994), Banach Center Publ., vol.~36,
  Polish Acad. Sci., Warsaw, 1996, pp.~25--44.

\bibitem[Bri97]{Brion97}
\bysame, \emph{Equivariant {C}how groups for torus actions}, Transform. Groups
  \textbf{2} (1997), no.~3, 225--267.

\bibitem[EG98a]{EdidinGraham98}
D.~Edidin and W.~Graham, \emph{Equivariant intersection theory}, Invent. Math.
  \textbf{131} (1998), no.~3, 595--634.

\bibitem[EG98b]{EdidinGraham98b}
\bysame, \emph{Localization in equivariant intersection theory and the {B}ott
  residue formula}, Amer. J. Math. \textbf{120} (1998), no.~3, 619--636.

\bibitem[FS97]{FultonSturmfels97}
W.~Fulton and B.~Sturmfels, \emph{Intersection theory on toric varieties},
  Topology \textbf{36} (1997), no.~2, 335--353.

\bibitem[Ful93]{Fulton93}
W.~Fulton, \emph{Introduction to toric varieties}, Annals of Mathematics
  Studies, vol. 131, Princeton University Press, Princeton, NJ, 1993.

\bibitem[Har92]{Harris92}
J.~Harris, \emph{Algebraic geometry}, Graduate Texts in Mathematics, vol. 133,
  Springer-Verlag, New York, 1992, A first course.

\bibitem[Kat07]{Katz06}
E.~Katz, \emph{A tropical toolkit}, preprint, math.AG/0610878v2, 2007.

\bibitem[Kim92]{Kimura92}
S.~Kimura, \emph{Fractional intersection and bivariant theory}, Comm. Algebra
  \textbf{20} (1992), no.~1, 285--302.

\bibitem[Kly90]{Klyachko90}
A.~Klyachko, \emph{Equivariant vector bundles on toral varieties}, Math.
  USSR-Izv. \textbf{35} (1990), no.~2, 337--375.

\bibitem[KM05]{KnutsonMiller05}
A.~Knutson and E.~Miller, \emph{Gr\"obner geometry of {S}chubert polynomials},
  Ann. of Math. (2) \textbf{161} (2005), no.~3, 1245--1318.

\bibitem[Law88]{Lawrence88}
J.~Lawrence, \emph{Valuations and polarity}, Discrete Comput. Geom. \textbf{3}
  (1988), no.~4, 307--324.

\bibitem[McM79]{McMullen78}
P.~McMullen, \emph{Lattice invariant valuations on rational polytopes}, Arch.
  Math. (Basel) \textbf{31} (1978/79), no.~5, 509--516.

\bibitem[Mik06]{Mikhalkin06}
G.~Mikhalkin, \emph{Tropical geometry and its applications}, preprint,
  math.AG/0601041, 2006.

\bibitem[MS05]{MillerSturmfels05}
E.~Miller and B.~Sturmfels, \emph{Combinatorial commutative algebra}, Graduate
  Texts in Mathematics, vol. 227, Springer-Verlag, New York, 2005.

\bibitem[Pay06a]{chowcohom}
S.~Payne, \emph{Equivariant {C}how cohomology of toric varieties}, Math. Res.
  Lett. \textbf{13} (2006), no.~1, 29--41.

\bibitem[Pay06b]{branchedcovers}
\bysame, \emph{Toric vector bundles, branched covers of fans, and the
  resolution property}, To appear in J. Alg. Geom. math.AG/0605537, 2006.

\bibitem[Pay07]{ehrhartseries}
\bysame, \emph{Ehrhart series and lattice triangulations}, To appear in Discr.
  Comput. Geom. math.CO/0702052, 2007.

\bibitem[Ros89]{Rossmann89}
W.~Rossmann, \emph{Equivariant multiplicities on complex varieties},
  Ast\'erisque (1989), no.~173-174, 11, 313--330, Orbites unipotentes et
  repr\'esentations, III.

\bibitem[Sam53]{Samuel53}
P.~Samuel, \emph{Alg\`ebre locale}, M\'emor. Sci. Math., no. 123,
  Gauthier-Villars, Paris, 1953.

\bibitem[Sam55]{Samuel55}
\bysame, \emph{M\'ethodes d'alg\`ebre abstraite en g\'eom\'etrie alg\'ebrique},
  Ergebnisse der Mathematik und ihrer Grenzgebiete (N.F.), Heft 4,
  Springer-Verlag, Berlin, 1955.

\end{thebibliography}
\bibliographystyle{amsalpha}

\end{document}